\documentclass[12pt]{article}

\usepackage{hyperref}

\usepackage{amsfonts, amsmath, amssymb, amsthm}
\usepackage{a4}
\usepackage{graphicx}

\newtheorem{theorem}{Theorem}

\newtheorem{proposition}{Proposition}

\theoremstyle{remark}
\newtheorem*{remark}{Remark}

\title{Heptagon relation in a direct sum}
\author{Igor G. Korepanov}

\date{}

\begin{document}

\sloppy

\maketitle

\begin{abstract}
An ansatz is proposed for heptagon relation, that is, algebraic imitation of five-dimen\-sional Pachner move 4--3. Our relation is realized in terms of matrices acting in a direct sum of one-dimen\-sional linear spaces corresponding to 4-faces.
\end{abstract}

\section{Introduction}\label{s:i}

Heptagon relation imitates algebraically a Pachner move 4--3 in a triangulation of a five-dimen\-sional piecewise linear (PL) manifold. This means a local re-building of the triangulation that leaves the manifold unchanged; to be more exact, it takes a cluster of four 5-simplices~$\Delta^5$ that form the \emph{star} of a 2-simplex (that is, a triangle~$\Delta^2$) and replaces it with a cluster of three 5-simplices that form the star of a 3-simplex (tetrahedron~$\Delta^3$). The fundamental Pachner theorem~\cite{Pachner,Lickorish} states (in particular) that, for a closed five-dimen\-sional PL manifold, any its triangulation can be transformed into any other one by a sequence of Pachner moves 4--3, 5--2, 6--1 and their inverses; the numbers here are of course numbers of 5-simplices before and after the move. In this paper, we restrict ourself to considering move 4--3.

Our algebraic imitation uses \emph{colorings} of \emph{four-dimen\-sional} faces (that is, \emph{pentachora}~$\Delta^4$) in our clusters: each~$\Delta^4$ is assigned a color which is an element of a given \emph{set of colors}~$X$. For each separate~$\Delta^5$, we define in Section~\ref{s:hg} a subset of \emph{permitted colorings} of its six 4-faces. For a cluster of 5-simplices, a coloring is by definition permitted if its restrictions onto all 5-simplices are permitted.

For any Pachner move, its initial and final clusters of simplices have the same boundary. For our move 4--3, it consists of twelve 4-simplices. By definition, heptagon relation holds if the sets of permitted \emph{boundary colorings}---the restrictions of permitted colorings onto the boundary---are the same for the two clusters. Such version of heptagon may be called \emph{set-theoretic}; the `direct-sum', or `matrix', relation introduced below is a particular case of this.

The contents of the remaining sections is as follows:
\begin{itemize}\itemsep 0pt
 \item in Section~\ref{s:hg}, we explain our `matrix' form of heptagon relation,
 \item in Section~\ref{s:ha}, we present our ansatz. That is, we write out algebraic formulas for matrix entries, leaving the explanation of their origin for the two next sections,
 \item in Section~\ref{s:he}, we introduce `edge vectors'---key algebraic structure behind our ansatz. Then we report the results of numerical experiments concerning these vectors (starting with the experimental fact that they exist),
 \item in Section~\ref{s:hc}, we construct edge vectors algebraically and prove that their properties agree with the experimental findings, and explain how our ansatz has been obtained,
 \item finally, in the concluding Section~\ref{s:d} we briefly discuss the obtained relations, including their connections with ``quantum'' relations involving tensor products of vector spaces, and possible generalizations.
\end{itemize}

\section{Direct-sum heptagon: generalities}\label{s:hg}

In our case, the set of colors will be, by definition, a \emph{field}~$X=F$. Colorings of, say, $n$~pentachora taken together form then a \emph{direct sum}~$F^n$ of $n$ copies of~$F$.

For the six faces of a separate 5-simplex, all colorings form the six-dimen\-sional linear space~$F^6$, and \emph{permitted colorings} are determined, by definition, by three \emph{linear relations} between the six colors. These relations are supposed to be generic enough, so that we can consider any three faces as `input' where we can assign any three colors; then the colors of the three remaining `output' faces are determined as their linear functions. We find it convenient to write both the `input' and `output' colors as 3-\emph{rows}, and specify the dependence between them by a $3\times 3$ matrix, acting on a row, of course, from the \emph{right}.

In writing out the heptagon relation, we can regard \emph{six} of the twelve boundary faces mentioned in Section~\ref{s:i} as `input' and six others as `output'; we write both input and output colors as 6-rows. Denote the $3\times 3$ matrix associated with one chosen 5-simplex as~$A$. Although~$A$ acts on only three colors, we can extend its action from 3-rows onto 6-rows in an obvious way: take its direct sum $A\oplus \mathbf 1$ with the identity matrix~$\mathbf 1$ acting on the three remaining colors. We will use the following notations: let $i,j,k=1\ldots 6$ be the three positions in the 6-rows on which matrix~$A$ actually acts, in this situation we denote $A\oplus \mathbf 1$ as~$A_{ijk}$. 

\begin{figure}
 \begin{center}
  \includegraphics[scale=1.25]{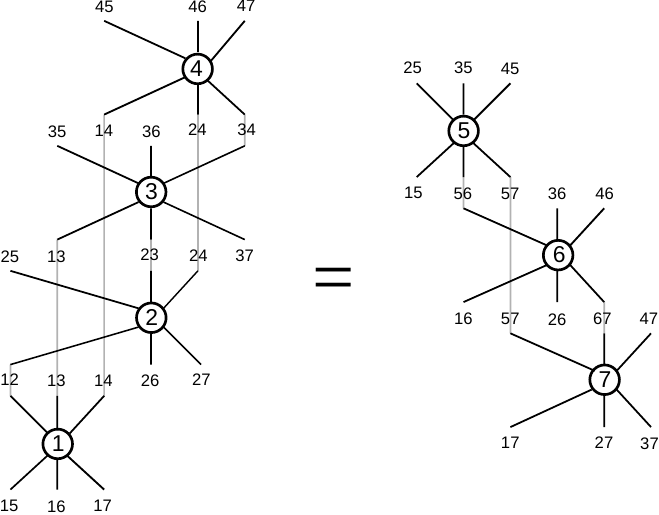}
 \end{center}
 \caption{Heptagon}
 \label{fig:h}
\end{figure}

Our heptagon relation involves, of course, seven such matrices $A=A^{(1)},\allowbreak \ldots,\allowbreak A^{(7)}$, and is as follows:
\begin{equation}\label{hepta}
A_{123}^{(1)}A_{145}^{(2)}A_{246}^{(3)}A_{356}^{(4)}=A_{356}^{(7)}A_{245}^{(6)}A_{123}^{(5)}.
\end{equation}
An explanation of why \eqref{hepta} is a right form for the heptagon can be deduced from comparing it to Equation~\cite[(4.7)]{DMH}. Alternatively, this can be explained using Figure~\ref{fig:h}.

Namely, in Figure~\ref{fig:h}, circles containing numbers $1,\ldots, 7$ represent the seven 5-simplices; the numbers are the same as parenthesized superscripts in~\eqref{hepta}. The l.h.s.\ of the figure represents the cluster of four 5-simplices, while the r.h.s.---the cluster of three. Edges represent \emph{4-faces}. If an edge coming from a circle is marked by two digits~$ij$, it means that this edge represents the 4-face common for 5-simplices $i$ and~$j$. This 4-face may be either inner for one cluster (as is, for instance, face~12), or be boundary---then it belongs to both clusters (like, for instance, face~15).

The `input' 4-faces correspond, in Figure~\ref{fig:h}, to the lower legs of a circle, while `output'---to the upper legs. The common boundary of the l.h.s.\ and r.h.s.\ of the Pachner move corresponds to the six `input' legs and six `output' legs of either l.h.s.\ or r.h.s.\ of Figure~\ref{fig:h}.

It can be checked that Figure~\ref{fig:h} represents indeed the 4--3 move as described in Section~\ref{s:i}. Recall that the subscripts $1,\ldots, 6$ in~\eqref{hepta} correspond simply to the position, counted from left to right, of a given 4-face in (the l.h.s.\ or r.h.s.\ of) Figure~\ref{fig:h}.

\section{The ansatz}\label{s:ha}

Below in Subsection~\ref{ss:ha}, we present our ansatz for heptagon relation. For its origins, and related algebraic structures, see Sections \ref{s:he} and~\ref{s:hc}.

Also, even before explaining these origins, we want to see whether our ansatz is specific just for the heptagon, or similar formulas can work for its analogues in other dimensions. Namely, we show in Subsection~\ref{ss:pa} that a similar ansatz also works at least for the three-dimen\-sional analogue of heptagon, that is, \emph{pentagon} relation.

\subsection{Explicit expressions for matrix entries}\label{ss:ha}

Let $F$ be a big enough field. We begin with introducing a triple $(\alpha_i,\beta_i,\gamma_i)$ of generic numbers $\alpha_i,\beta_i,\gamma_i\in F$ for each~$i=1,\ldots,7$. Then, we introduce determinants
\begin{equation}\label{dijk}
d_{ijk} = \left| \begin{matrix} \alpha_i & \alpha_j & \alpha_k \\ \beta_i & \beta_j & \beta_k \\ \gamma_i & \gamma_j & \gamma_k \end{matrix} \right|.
\end{equation}
We assume that our alphas, betas and gammas are \emph{generic} enough: to be exact, this will mean that $d_{ijk}$ must not vanish for any pairwise different $i,j,k$.

\begin{figure}
 \begin{center}
  \includegraphics[scale=1.25]{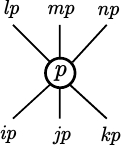}
 \end{center}
 \caption{Matrix~$A^{(p)}$ for heptagon}
 \label{fig:Ahepta}
\end{figure}

Our ansatz consists in setting the entry of matrix~$A^{(p)}$ corresponding to the input (lower) leg $ip=pi$ and output (upper) leg $lp=pl$ to be
\begin{equation}\label{ansatz-h}
\left( A^{(p)} \right)_{ip}^{lp} = \frac{d_{jlp}d_{klp}}{d_{ijp}d_{ikp}},
\end{equation}
where $jp$ and~$kp$ are other input legs of~$A^{(p)}$. That is, $A^{(p)}$ can look as in Figure~\ref{fig:Ahepta},
or be obtained from that by permutations of lower legs and upper legs separately.

For instance, this means that
\[
A^{(7)} = \begin{pmatrix} 
 \dfrac{d_{257}d_{357}}{d_{127}d_{137}} & \dfrac{d_{267}d_{367}}{d_{127}d_{137}} 
               & \dfrac{d_{247}d_{347}}{d_{127}d_{137}} \\[2.5ex]
 \dfrac{d_{157}d_{357}}{d_{217}d_{237}} & \dfrac{d_{167}d_{367}}{d_{217}d_{237}} 
               & \dfrac{d_{147}d_{347}}{d_{217}d_{237}} \\[2.5ex]
 \dfrac{d_{157}d_{257}}{d_{317}d_{327}} & \dfrac{d_{167}d_{267}}{d_{317}d_{327}} 
               & \dfrac{d_{147}d_{247}}{d_{317}d_{327}}
\end{pmatrix} .
\]

\begin{remark}
Don't forget that our matrices act on \emph{rows}, that is, from the right!
\end{remark}

Given~\eqref{ansatz-h}, heptagon relation~\eqref{hepta} can be proved by brute force: a direct calculation using computer algebra. Not being content with just this, we provide below a more conceptual proof, based on Pl\"ucker bilinear relations~\cite[p.~211]{Pluecker}.

\subsection{Pentagon, or a lite version of ansatz}\label{ss:pa}

\begin{figure}[p]
\begin{center}
\includegraphics[scale=1.25]{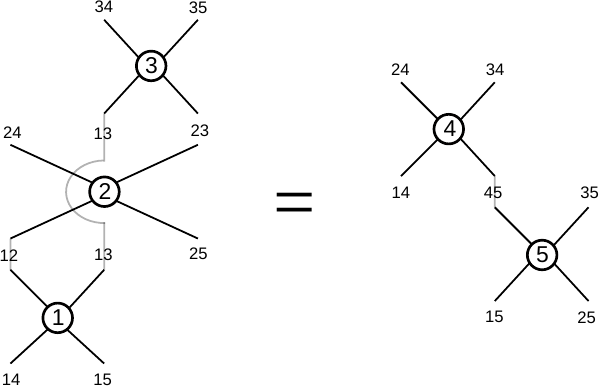}
\end{center}
\caption{Pentagon}
\label{fig:p}
\end{figure}

\begin{figure}[p]
 \begin{center}
  \includegraphics[scale=1.25]{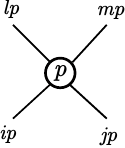}
 \end{center}
 \caption{Matrix~$A^{(p)}$ for pentagon}
 \label{fig:Apenta}
\end{figure}

Pentagon relation, that is, a three-dimen\-sional analogue of heptagon~\eqref{hepta}, looks as follows:
\begin{equation}\label{penta}
A_{12}^{(1)} A_{13}^{(2)} A_{23}^{(3)} = A_{23}^{(5)} A_{12}^{(4)} .
\end{equation}
A graphic representation of this can be seen in Figure~\ref{fig:p}.
And the analogue of ansatz~\eqref{ansatz-h} is now simply as follows:
\begin{equation}\label{ansatz-p}
\left( A^{(p)} \right)_{ip}^{lp} = \frac{d_{jlp}}{d_{ijp}},
\end{equation}
with the same $d_{ijk}$~\eqref{dijk} as before, and $A^{(p)}$ as in Figure~\ref{fig:Apenta}.
Again, the validity of~\eqref{penta} can be checked directly.

\begin{remark}
Actually, ansatz~\eqref{ansatz-p} is a simple variation on the theme of the pentagon relations found in~\cite{pentagon}, but we will not go into those details here.
\end{remark}

\section{Edge vectors and the experiment}\label{s:he}

The key algebraic structures that led us to ansatz~\eqref{ansatz-h} were \emph{edge vectors}.

\subsection{Definition of edge vectors}\label{ss:ed}

By edge vector~$e_b$, corresponding to an edge
\begin{equation}\label{bijji}
b=ij=ji
\end{equation}
joining vertices $i$ and~$j$, we mean here a \emph{nonzero permitted coloring} of both the l.h.s.\ (\,=\,initial cluster of four 5-simplices) and r.h.s.\ (\,=\,resulting cluster of three 5-simplices) of move 4--3 satisfying the following conditions:
\begin{itemize}\itemsep 0pt
 \item[(a)] only those pentachora~$u$ that contain edge~$b$ may be colored by a nonzero element of field~$F$,
 \item[(b)] the restrictions of the l.h.s.\ and r.h.s\ colorings onto their common boundary are the same.
\end{itemize}

As the l.h.s.\ and r.h.s\ of an $n$-dimen\-sional Pachner move are known to form together the boundary~$\partial \Delta^{n+1}$ of a simplex of the next dimension, $\partial \Delta^6$ in our case, edge vector can be thought of as a coloring of~$\partial \Delta^6$ satisfying~(a). We will use this point of view when it is convenient.

The reader may have noticed also that we denote 6-vertex simplices (that is, \emph{5-simplices}~$\Delta^5$\,!) by the same letters $i,j,\ldots$ as vertices. This will be convenient for us as long as we work within (the boundary of) one 6-simplex with vertices $1,\ldots,7$, that is, $\Delta^6=1234567$. Specifically, $i=1,\ldots,7$ denotes below either vertex~$i$, or the 6-vertex simplex containing all vertices \emph{except}~$i$ (that is, `1' may denote $234567$. The exact meaning will hopefully be clear from the context).

\subsection{Example}\label{ss:ee}

For illustration of how an edge vector may look in terms of Figure~\ref{fig:h}, we take~$e_{47}$ as an example.

According to what we said above about our notations for vertices and their complementary 6-vertex simplices, 4-faces containing edge~47 are exactly those corresponding to edges in Figure~\ref{fig:h} marked by two digits \emph{neither} of which is 4 or~7. So, wherever 4 or~7 \emph{is} present at a line, there must be a zero.

\begin{figure}
 \begin{center}
  \includegraphics[scale=1.25]{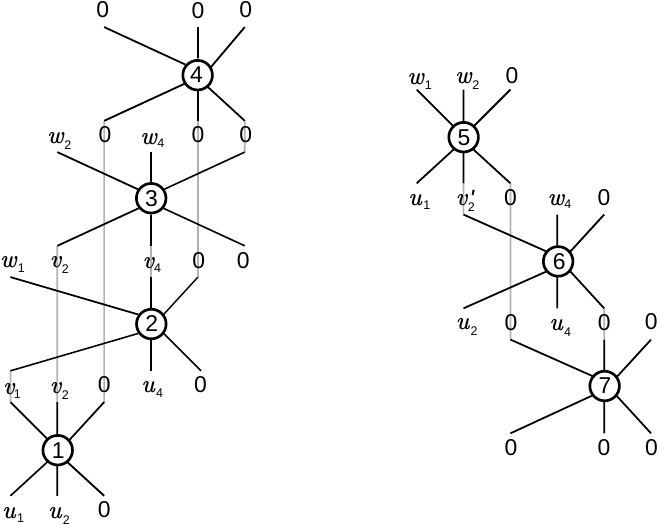}
 \end{center}
 \caption{Components of edge vector $e_{47}$}
 \label{fig:l}
\end{figure}

The components of~$e_{47}$ are depicted at the edges in Figure~\ref{fig:l}.
We see that the existence of~$e_{47}$ imposes some \emph{compatibility} conditions on matrices~$A^{(p)}$. For instance, we see in the left-hand side of Figure~\ref{fig:l} that there are such numbers $u_1,\ldots,w_3 \in F$ (not all zero) that
\begin{align}
\begin{pmatrix} u_1 & u_2 & 0 \end{pmatrix}A^{(1)} &= \begin{pmatrix} v_1 & v_2 & 0 \end{pmatrix}, \label{l1} \\
\begin{pmatrix} v_1 & u_4 & 0 \end{pmatrix}A^{(2)} &= \begin{pmatrix} w_1 & v_4 & 0 \end{pmatrix}, \label{l2} \\
\begin{pmatrix} v_2 & v_4 & 0 \end{pmatrix}A^{(3)} &= \begin{pmatrix} w_2 & w_4 & 0 \end{pmatrix}, \label{l3}
\end{align}
and relations \eqref{l1}--\eqref{l3} can be interpreted as some sort of compatibility between matrices $A^{(1)}$, $A^{(2)}$ and~$A^{(3)}$, 

\subsection{Experiment}\label{ss:xpr}

Relation~\eqref{hepta} is a system of nonlinear equations on the entries of seven matrices~$A_{ijk}^{(p)}$. There are efficient algorithms for solving such systems numerically, at least for the field $F=\mathbb R$ of real numbers.

Specifically, we used the \emph{Levenberg--Marquardt} algorithm~\cite{Levenberg,Marquardt}, starting from randomly chosen initial values of matrix entries and arriving at a high precision solution of~\eqref{hepta}.

All such solutions turned out to have the following properties:
\begin{itemize}\itemsep 0pt
 \item[(i)] edge vectors for all edges of~$\partial \Delta^6$ do exist,
 \item[(ii)] choose a pair consisting of a 4-face~$u$ of $\partial\Delta^6$ and a vertex of~$u$. For instance, let $u=ijklm$, and the chosen vertex be~$i$. In this situation, there is a \emph{linear dependence} between $e_{ij}$, $e_{ik}$, $e_{il}$ and~$e_{im}$, that is, vectors belonging to edges~$b$ such that $i\in b\subset u$, while any three of these vectors are linearly independent. We write this linear dependence as follows:
\begin{equation}\label{lmnx}
\lambda_{i,ij}^{(u)} e_{ij}+\lambda_{i,ik}^{(u)} e_{ik}+\lambda_{i,il}^{(u)} e_{il}+\lambda_{i,im}^{(u)} e_{im}=0
\end{equation}
(so, the parenthesized superscript and the first subscript of a lambda mean the pentachoron and its vertex to which the linear relation belongs, while the second subscript means the edge),
 \item[(iii)] in the above notations, for a given $u=ijklm$, there is a linear dependence \emph{between the linear dependences} of type~\eqref{lmnx} in vertices $i,j,k,l,m$. This means that if we normalize (multiply by nonzero numbers) these dependences properly and add them all together, the coefficients of all~$e_b$ will vanish. As one can easily see, this vanishing looks as follows:
\begin{equation}\label{ll}
\lambda_{i,ij}^{(u)}+\lambda_{j,ij}^{(u)}=0, \qquad \ldots\,, \qquad \lambda_{l,lm}^{(u)}+\lambda_{m,lm}^{(u)}=0,
\end{equation}
 \item[(iv)] when restricted to one $\Delta^5$, edge vectors generate the whole 3-dimen\-sional space of its permitted colorings,
 \item[(v)] when taken as they are, edge vectors generate the whole 6-dimen\-sional space of permitted colorings of~$\partial \Delta^6$.
\end{itemize}

\section{Construction of edge vectors}\label{s:hc}

Our task now is to construct algebraically---not just numerically---edge vectors satisfying items (i)--(v) of Subsection~\ref{ss:xpr}, and then the corresponding heptagon relations. The logical structure of this section will be as follows: 
\begin{itemize}\itemsep 0pt
 \item we \emph{assume} that items (i)--(v) hold, and moreover make one more assumption---see~\eqref{zvezda} below---and derive formulas for coefficients~$\lambda$ of linear dependences~\eqref{lmnx}. Happily, these linear dependences turn out to be normalizable in such way that~\eqref{ll} holds. This is done in Subsection~\ref{ss:l},
 \item then, in Subsection~\ref{ss:e}, we derive an explicit expression for the components of edge vectors,
 \item our explicit expressions followed \emph{necessarily} from (i)--(v) and~\eqref{zvezda} (and our normalization choices for edge vectors and $\lambda$'s, to be pedantic). But now we can check that the edge vectors and $\lambda$'s given by the obtained expressions satisfy (i)--(v) indeed; this statement is clarified and partially proved in Subsection~\ref{ss:lb},
 \item and finally, in Subsection~\ref{ss:r}, we finish the mentioned proof, and derive the heptagon relation and formula~\eqref{ansatz-h} for matrix entries.
\end{itemize}

\subsection{Coefficients of linear dependences between edge vectors}\label{ss:l}

Linear dependences~\eqref{lmnx}, together with the fact that any three of the vectors in~\eqref{lmnx} are linearly independent, mean that the six edge vectors~$e_{ij}$, taken for a given vertex~$i$ and all $j\ne i$, \ $1\le j\le 7$, span a 3-dimen\-sional linear space~$V_i$. Take a basis in~$V_i$; we would like to denote its three vectors as $e_i^{(\alpha)}$, $e_i^{(\beta)}$ and~$e_i^{(\gamma)}$. Each~$e_{ij}$ is a linear combination of these, and a \emph{minimalistic}---and productive!---idea is to assume that its coefficients depend \emph{only on}~$j$ and not on~$i$:
\begin{equation}\label{zvezda}
e_{ij} = \alpha_j e_i^{(\alpha)} + \beta_j e_i^{(\beta)} + \gamma_j e_i^{(\gamma)}.
\end{equation}

Numbers $\alpha_i,\beta_i,\gamma_i$ in~\eqref{zvezda} will, of course, turn out very soon to be the same as those introduced in the beginning of Subsection~\ref{ss:ha} (don't forget that indices $i$ or~$j$ can denote both a vertex and its complementary 6-vertex simplex, see the last paragraph of Subsection~\ref{ss:ed}).

\begin{proposition}
Suppose that edge vectors for our heptagon relation exist and are given by~\eqref{zvezda}, with $\alpha_i,\beta_i,\gamma_i$ generic in the same exact sense as stated after~\eqref{dijk}. Suppose also that there are linear dependences~\eqref{lmnx}, and no more linear dependences between the same vectors. Then, coefficients~$\lambda$ can be chosen so that they also obey~\eqref{ll}, and this determines them uniquely up to one overall factor for each pentachoron~$u$. With a natural choice of these factors, the lambdas are as follows:
\begin{equation}\label{ijij}
\lambda_{i,ij}^{(u)} = - \lambda_{j,ij}^{(u)} = \epsilon_{ijklm} \left| \begin{matrix} \alpha_k & \alpha_l & \alpha_m \\ \beta_k & \beta_l & \beta_m \\ \gamma_k & \gamma_l & \gamma_m \end{matrix} \right| = \epsilon_{ijklm} d_{klm},
\end{equation}
where $i,j,k,l,m$ are the vertices of~$u$ (they are thus pairwise different integers between 1 and~7, but no assumption is made about the order in which these numbers are taken), and $\epsilon_{ijklm}$ is the sign
\[
\epsilon_{ijklm}=(-1)^{\mathrm{parity}}=\pm 1
\]
of the permutation between numbers $i,j,k,l,m$ taken in this order and the increasing order ($d_{klm}$ is of course defined according to~\eqref{dijk}).
\end{proposition}

\begin{proof}
For a given vertex~$i$ and pentachoron~$u\ni i$, linear dependence~\eqref{lmnx} becomes a system of three linear homogeneous equations for the lambdas, whose solution is determined up to an overall factor. Fix pentachoron~$u$, then there is one such factor for each of its five vertices $i\in u$; writing out the corresponding solutions, we arrive at once at the proposition statement.
\end{proof}

\subsection{Expression for edge vectors}\label{ss:e}

Next, we are going to derive explicit formulas for the components~$e_b|_u$. Here $b$ is an edge, $u$ is a 4-face (pentachoron), and $e_b|_u$ is the color in which edge vector~$e_b$ paints~$u$. First, we show how to calculate the \emph{ratio} of the components $e_{b_1}|_u$ and~$e_{b_2}|_u$ corresponding to two edges $b_1\subset u$ and $b_2\subset u$ having one common vertex~$i=1$. For clarity, we do it on a typical example, namely, let $u=12345$, \ $b_1=12$, \ $b_2=13$, so $i=1$.

Write out linear dependence~\eqref{lmnx} for a \emph{different} pentachoron, namely $u'=ijklm=12367$, that is, containing all three vertices of our chosen edges $b_1$ and~$b_2$, but no other vertices of~$u$. Taking into account that
\[
e_{16}|_{12345}=e_{17}|_{12345}=0,
\]
because neither edge~$16$ nor~$17$ belongs to~$u$, we get
\begin{equation}\label{l/l}
\lambda_{1,12}e_{12}|_{12345}+\lambda_{1,13}e_{13}|_{12345}=0.
\end{equation}
Wherefrom it follows, on substituting lambdas from~\eqref{ijij}, that
\begin{equation}\label{d/d}
\frac{e_{12}|_{12345}}{e_{13}|_{12345}} = \frac{d_{267}}{d_{367}}.
\end{equation}

Equality~\eqref{d/d}---being a typical example---admits any permutation of numbers $1,\ldots,7$. Taking also into account the symmetry
\[
e_{ij}=e_{ji}
\]
(recall~\eqref{bijji}), we see that, for a given~$u$, all ratios $\dfrac{e_{b_1}|_u}{e_{b_2}|_u}$ for edges $b_1,b_2\subset u$ are fixed, so, all $e_b|_u$ are fixed up to one possible common factor depending only on~$u$. It may be called \emph{gauge} factor; it is simply responsible for a possible re-scaling of the colors on face~$u$, which does not affect the heptagon relation.

\begin{proposition}
With a proper choice of gauge factors, the explicit expression for components~$e_b|_u$ is as follows:
\begin{equation}\label{dd}
e_{ij}|_u = d_{ilm}d_{jlm},
\end{equation}
where $l$ and~$m$ are the two vertices \emph{not} belonging to~$u$ (while $i$ and~$j$, of course, do belong).
\end{proposition}

\begin{proof}
It remains to note that \eqref{dd} is compatible with all relations of type~\eqref{d/d}:
\end{proof}

Note, by the way, that the condition~(a) of Subsection~\ref{ss:ed}, that is,
\[
e_b|_u=0\quad \text{for}\quad b\not\subset u,
\]
follows from~\eqref{dd} automatically!

\subsection[Existence: checking (i)--(v) for the obtained edge vectors and coefficients~$\lambda$]{Existence: checking (i)--(v) for the obtained edge vectors and coefficients~$\boldsymbol{\lambda}$}\label{ss:lb}

We now take~\eqref{dd} as the \emph{definition} of vectors~$e_b$, and~\eqref{ijij} as the definition of lambdas, and check all the items (i)--(v) of Subsection~\ref{ss:xpr}, in the sense that we are going to explain.

First, a clarification must be made about item~(i): it meant that edge vectors exist \emph{for a heptagon relation} obtained experimentally. Now that we have explicit vectors~\eqref{dd}, clearly satisfying the two conditions in Subsection~\ref{ss:ed}, what must be shown for them is that they really \emph{yield the heptagon relation}. We will do this later, in Subsection~\ref{ss:r}, while here we are going to check (ii), (iii), (iv), and partially~(v); see the following propositions for the exact sense of this.

\begin{proposition}[Checking item~(ii)]
For the constructed edge vectors, linear dependences~\eqref{lmnx} hold, and there are no more linear dependences between the same four vectors.
\end{proposition}

\begin{proof}
Write out the component of the l.h.s.\ of~\eqref{lmnx} corresponding to a face~$u'$ (that may coincide or not with~$u$ in~\eqref{lmnx}). After factoring out the common factor~$d_{il'm'}$, where $l'$ and~$m'$ are the two vertices \emph{not belonging} to~$u'$, one is left with nothing but the well-known \emph{Pl\"ucker bilinear relation}~\cite[p.~211]{Pluecker}.

That there is no more linear dependence between the same vectors, follows simply from the fact that equalities of type~\eqref{l/l} fix unambiguously all ratios between the relevant coefficients~$\lambda$.
\end{proof}

\begin{proposition}[Checking item~(iii)] 
For the constructed lambdas, equalities~\eqref{ll} hold.
\end{proposition}

\begin{proof}
This is obvious.
\end{proof}

\begin{proposition}[Checking item~(iv)]\label{p:iv}
When restricted to one $\Delta^5$, the constructed edge vectors generate a 3-dimen\-sional linear space of its colorings.
\end{proposition}

And these colorings will be called, by definition, permitted colorings of~$\Delta^5$ induced by our constructed system of edge vectors.

\begin{proof}
We show that, for instance, restrictions of $e_{12}$, $e_{13}$ and~$e_{23}$ onto the faces of 5-simplex~$123456$---denoted below as $e_{12}|_{123456}$, etc.---form a basis in the space of all colorings that arise as restrictions onto~$123456$ of any linear combinations of \emph{all} edge vectors.

First, suppose there is a linear relation
\begin{equation}\label{nu}
(\nu _{12}e_{12}+\nu _{13}e_{13}+\nu _{23}e_{23})|_{123456} = 0.
\end{equation}
All $e_{ij}$ in~\eqref{nu} except~$e_{12}$ have zero components corresponding to face~$12456$, which implies immediately $\nu _{12} =0$, and similarly we get $\nu _{13} =0$ and $\nu _{23} =0$.

On the other hand, note that if we put $m=7$ in~\eqref{lmnx}, then the last term in the restriction of~\eqref{lmnx} onto $123456$ vanishes. So, what remains in the mentioned restriction allows to express any of the remaining~$e_{ij}|_{123456}$ through those~$e_{ij}|_{123456}$ in~\eqref{nu}, using this restriction of~\eqref{lmnx} with properly chosen $i,j,k,l$, and $m=7$, one or several times. 
\end{proof}

\begin{proposition}[Checking item~(v) partially]\label{p:v}
The constructed edge vectors generate a 6-dimen\-sional space of permitted colorings of~$\partial \Delta^6$.
\end{proposition}

Recall that a permitted coloring of~$\partial \Delta^6$ is defined as such a coloring whose restriction on any~$\Delta^5$ is permitted. So, the space generated by edge vectors consists indeed of permitted colorings of~$\partial \Delta^6$; what is \emph{not} (yet) stated in Proposition~\ref{p:v} is that \emph{all} permitted colorings of~$\partial \Delta^6$ are obtained that way. What \emph{does} follow from this proposition is that the space of permitted colorings is \emph{at least} 6-dimen\-sional.

\begin{proof}
  Six linearly independent edge vectors are, for instance, $e_{12}$, $e_{13}$, $e_{14}$, $e_{23}$, $e_{24}$, and~$e_{34}$---that is, corresponding to edges drawn between four chosen vertices, $1,\ldots,4$ in our case.

Indeed, suppose there is a linear relation
\begin{equation}\label{mu}
\mu _{12}e_{12}+\mu _{13}e_{13}+\mu _{14}e_{14}+\mu _{23}e_{23}+\mu _{24}e_{24}+\mu _{34}e_{34} = 0.
\end{equation}
All $e_{ij}$ in~\eqref{mu} except~$e_{12}$ have zero components corresponding to face~$12567$, which implies immediately $\mu _{12} =0$, and the same can obviously be done for all~$\mu$'s.

On the other hand, any of the remaining~$e_{ij}$ can be expressed through those in~\eqref{mu}, using~\eqref{lmnx} with properly chosen $i,j,k,l,m$ one or several times. 
\end{proof}

\subsection{Explicit form of matrix entries, and the heptagon relation}\label{ss:r}

Consider once again Figure~\ref{fig:Ahepta}. It corresponds to the 5-simplex~$ijklmn$, which we also characterize alternatively as 5-simplex \emph{without} vertex~$p$; here $i,j,k,l,m,n,p$ make, of course, a permutation of numbers $1,2,3,4,5,6,7$. The legs correspond to 4-faces also without the two corresponding vertices.

Consider now edge vector~$e_{jk}$. As 4-faces denoted $jp$ and~$kp$ do \emph{not} contain $j$ and~$k$, respectively, the corresponding components of~$e_{jk}$ vanish. Other components are given by formula~\eqref{dd}, so, the components of~$e_{jk}$ on the three lower legs of Figure~\ref{fig:Ahepta}---call them together, for a moment, $e_{jk}|_{\mathrm{lower}}$---are
\begin{equation}\label{ejk-niz}
e_{jk}|_{\mathrm{lower}} = \begin{pmatrix} d_{jip}d_{kip} & 0 & 0 \end{pmatrix},
\end{equation}
while the components on the three upper legs are
\begin{equation}\label{ejk-verh}
e_{jk}|_{\mathrm{upper}} = \begin{pmatrix} d_{jlp}d_{klp} & d_{jmp}d_{kmp} & d_{jnp}d_{knp} \end{pmatrix}.
\end{equation}

Similar expressions hold also for $e_{ik}$ and~$e_{ij}$, and we see once again (compare the proof of Proposition~\ref{p:iv}) that the restrictions of $e_{jk}$, $e_{ik}$ and~$e_{ij}$ on our 5-simplex~$p$ form a basis in its permitted colorings. As we know from Proposition~\ref{p:iv} that they form exactly a 3-dimen\-sional space, this space can indeed be described as follows: take any colors for the three lower legs, arrange them in a row, then the colors on the upper legs are that row multiplied on the right by~$A^{(p)}$, where $A^{(p)}$ is determined from
\begin{equation}\label{Ap}
\begin{pmatrix} e_{jk}|_{\mathrm{lower}} \\ e_{ik}|_{\mathrm{lower}} \\ e_{ij}|_{\mathrm{lower}} \end{pmatrix} A^{(p)} = \begin{pmatrix} e_{jk}|_{\mathrm{upper}} \\ e_{ik}|_{\mathrm{upper}} \\ e_{ij}|_{\mathrm{upper}} \end{pmatrix} .
\end{equation}
Note the diagonal form of the leftmost matrix in~\eqref{Ap}. This allows to calculate~$A^{(p)}$ at once, and we have arrived at the following proposition.

\begin{proposition}
The permitted colorings of a~$\Delta^5$ induced by the constructed edge vectors can be described in terms of matrix~$A^{(p)}$ with entries given by~\eqref{ansatz-h}.
 \qed
\end{proposition}

We can now describe the linear space of permitted colorings of the l.h.s.\ of move 4--3 as follows: assign any colors to the six ``lower input'' legs in the l.h.s.\ of Figure~\ref{fig:h}, and obtain the colors on other legs by consecutive action of $A^{(1)}$, $A^{(2)}$, $A^{(3)}$ and~$A^{(4)}$. This space is thus proved to be six-dimen\-sional. Similar description holds also for the r.h.s. The space of permitted colorings of~$\partial \Delta^6$---the union of the l.h.s.\ and r.h.s.---must be at least 6-dimensional, as Proposition~\ref{p:v} tells us, and this can be only in the case where heptagon relation~\eqref{hepta} holds. We have thus concluded the proof of the following theorem, without referring to computer algebra.

\begin{theorem}
Matrices~$A^{(p)}$ defined according to~\eqref{ansatz-h} satisfy heptagon relation~\eqref{hepta}.
 \qed
\end{theorem}

\section{Discussion}\label{s:d}

\subsection{Other types of similar relations}

Since the discovery of the heptagon relations described in this paper, one more type of such relations has been discovered~\cite{DK}. Whether these two types of heptagon relations can be included in a single, and more general, scheme, remains unknown at the moment when these lines are being written.

Also, it may be interesting to study possible links with similar---namely, \emph{hexagon}---relations parameterized by \emph{simplicial cocycles}~\cite{nonconstant}.

\subsection{Polygon and simplex relations}

Heptagon relation is one of the \emph{polygon relations}, and these are often studied together with \emph{simplex relations}~\cite{DMH}. There are many reasons; we just note here that one particular indication to some kinship between polygon and simplex relations is that the \emph{left-hand sides} of heptagon and \emph{tetrahedron} are virtually the same, compare~\cite[Eq.~(4.7)]{DMH} with the first of two unnumbered equations on~\cite[page~17]{DMH}.

It must be said, however, that there are different versions of both polygon and simplex relations (and different versions may be written in the same symbolic form!). We briefly described \emph{set-theoretic} heptagon in Section~\ref{s:i}, and then focused on its particular, and maybe most important, case---direct-sum relation. Most popular version of all these relations seems to be, however, their \emph{quantum}, or \emph{tensor}, version, where our direct sums of vector spaces are replaced with tensor products. In this connection, we would like to make here three following remarks:
\begin{itemize} \itemsep 0pt
 \item an interesting study of direct-sum \emph{simplex} relations has been done by Hietarinta~\cite{Hietarinta} in 1997,
 \item Hietarinta also explains how to make (simple) quantum relations from their direct-sum (or even set-theoretic) versions,
 \item more quantum relations can be obtained if we add \emph{cohomology} of direct-sum relations, like it was done in~\cite{cubic,nonconstant} for \emph{hexagon}.
\end{itemize}

\subsection{Possible generalizations}

Finally, expressions \eqref{ansatz-p} and~\eqref{ansatz-h} obviously suggest a generalization for higher $n$-gon relations with odd~$n$. This will be the subject of our future work. Right here, we can announce that the role of edges passes, for an arbitrary odd $n\ge 5$, to $\frac{n-5}{2}$-simplices, and to each such simplex~$\Delta^{(n-5)/2}$ belongs a vector (defined in obvious generalization of our edge vectors) whose components are products of determinants~\eqref{dijk} over the vertices of~$\Delta^{(n-5)/2}$:
\begin{equation*}
e_{\Delta^{(n-5)/2}}|_v = \prod_{i\in \Delta^{(n-5)/2}} d_{ilm}.
\end{equation*}
Here, $v$ is an $(n-3)$-face, $l$ and~$m$ are the two vertices \emph{not} belonging to~$v$, and the whole expression clearly generalizes~\eqref{dd}.

\subsection*{Acknowledgments}
\emph{Numeric calculations}: Alexey Korepanov taught me about the Levenberg--Marquardt algorithm, as well as how to use a relevant \texttt{C++} program.

\emph{Symbolic calculations} were made using wxMaxima.

\end{document}